\documentclass[twoside]{article}
\usepackage{theorem,amsmath,amssymb}

\makeatletter
\newcount\@hour\newcount\@minute\chardef\@x10\chardef\@xv60
\def\tcitime{
\def\@time{%
  \@minute\time\@hour\@minute\divide\@hour\@xv
  \ifnum\@hour<\@x 0\fi\the\@hour:%
  \multiply\@hour\@xv\advance\@minute-\@hour
  \ifnum\@minute<\@x 0\fi\the\@minute
  }}%

\@ifundefined{hyperref}{}{}

\@ifundefined{qExtProgCall}{\def\qExtProgCall#1#2#3#4#5#6{\relax}}{}
\def\QCTOpt[#1]#2{%
  \def\QCTOptB{#1}
  \def\QCTOptA{#2}
}
\def\QCTNOpt#1{%
  \def\QCTOptA{#1}
  \let\QCTOptB\empty
}
\def\Qct{%
  \@ifnextchar[{%
    \QCTOpt}{\QCTNOpt}
}
\def\QCBOpt[#1]#2{%
  \def\QCBOptB{#1}
  \def\QCBOptA{#2}
}
\def\QCBNOpt#1{%
  \def\QCBOptA{#1}
  \let\QCBOptB\empty
}
\def\Qcb{%
  \@ifnextchar[{%
    \QCBOpt}{\QCBNOpt}
}
\def\PrepCapArgs{%
  \ifx\QCBOptA\empty
    \ifx\QCTOptA\empty
      {}%
    \else
      \ifx\QCTOptB\empty
        {\QCTOptA}%
      \else
        [\QCTOptB]{\QCTOptA}%
      \fi
    \fi
  \else
    \ifx\QCBOptA\empty
      {}%
    \else
      \ifx\QCBOptB\empty
        {\QCBOptA}%
      \else
        [\QCBOptB]{\QCBOptA}%
      \fi
    \fi
  \fi
}
\newcount\GRAPHICSTYPE
\GRAPHICSTYPE=\z@
\def\GRAPHICSPS#1{%
 \ifcase\GRAPHICSTYPE
   \special{ps: #1}%
 \or
   \special{language "PS", include "#1"}%
 \fi
}%
\def\graffile#1#2#3#4{%
    \leavevmode
    \raise -#4 \BOXTHEFRAME{%
        \hbox to #2{\raise #3\hbox to #2{\null #1\hfil}}}%
}%
\def\draftbox#1#2#3#4{%
 \leavevmode\raise -#4 \hbox{%
  \frame{\rlap{\protect\tiny #1}\hbox to #2%
   {\vrule height#3 width\z@ depth\z@\hfil}%
  }%
 }%
}%
\newcount\draft
\draft=\z@

\newif\ifwasdraft
\wasdraftfalse

\def\GRAPHIC#1#2#3#4#5{%
 \ifnum\draft=\@ne\draftbox{#2}{#3}{#4}{#5}%
  \else\graffile{#1}{#3}{#4}{#5}%
  \fi
 }%
\def\addtoLaTeXparams#1{%
    \edef\LaTeXparams{\LaTeXparams #1}}%
\newif\ifBoxFrame \BoxFramefalse
\newif\ifOverFrame \OverFramefalse
\newif\ifUnderFrame \UnderFramefalse

\def\BOXTHEFRAME#1{%
   \hbox{%
      \ifBoxFrame
         \frame{#1}%
      \else
         {#1}%
      \fi
   }%
}

\def\doFRAMEparams#1{\BoxFramefalse\OverFramefalse\UnderFramefalse\readFRAMEparams#1\end}%
\def\readFRAMEparams#1{%
 \ifx#1\end%
  \let\next=\relax
  \else
  \ifx#1i\dispkind=\z@\fi
  \ifx#1d\dispkind=\@ne\fi
  \ifx#1f\dispkind=\tw@\fi
  \ifx#1t\addtoLaTeXparams{t}\fi
  \ifx#1b\addtoLaTeXparams{b}\fi
  \ifx#1p\addtoLaTeXparams{p}\fi
  \ifx#1h\addtoLaTeXparams{h}\fi
  \ifx#1X\BoxFrametrue\fi
  \ifx#1O\OverFrametrue\fi
  \ifx#1U\UnderFrametrue\fi
  \ifx#1w
    \ifnum\draft=1\wasdrafttrue\else\wasdraftfalse\fi
    \draft=\@ne
  \fi
  \let\next=\readFRAMEparams
  \fi
 \next
 }%
\def\IFRAME#1#2#3#4#5#6{%
      \bgroup
      \let\QCTOptA\empty
      \let\QCTOptB\empty
      \let\QCBOptA\empty
      \let\QCBOptB\empty
      #6%
      \parindent=0pt%
      \leftskip=0pt
      \rightskip=0pt
      \setbox0 = \hbox{\QCBOptA}%
      \@tempdima = #1\relax
      \ifOverFrame
          \typeout{This is not implemented yet}%
          \show\HELP
      \else
         \ifdim\wd0>\@tempdima
            \advance\@tempdima by \@tempdima
            \ifdim\wd0 >\@tempdima
               \textwidth=\@tempdima
               \setbox1 =\vbox{%
                  \noindent\hbox to \@tempdima{\hfill\GRAPHIC{#5}{#4}{#1}{#2}{#3}\hfill}\\%
                  \noindent\hbox to \@tempdima{\parbox[b]{\@tempdima}{\QCBOptA}}%
               }%
               \wd1=\@tempdima
            \else
               \textwidth=\wd0
               \setbox1 =\vbox{%
                 \noindent\hbox to \wd0{\hfill\GRAPHIC{#5}{#4}{#1}{#2}{#3}\hfill}\\%
                 \noindent\hbox{\QCBOptA}%
               }%
               \wd1=\wd0
            \fi
         \else
            \ifdim\wd0>0pt
              \hsize=\@tempdima
              \setbox1 =\vbox{%
                \unskip\GRAPHIC{#5}{#4}{#1}{#2}{0pt}%
                \break
                \unskip\hbox to \@tempdima{\hfill \QCBOptA\hfill}%
              }%
              \wd1=\@tempdima
           \else
              \hsize=\@tempdima
              \setbox1 =\vbox{%
                \unskip\GRAPHIC{#5}{#4}{#1}{#2}{0pt}%
              }%
              \wd1=\@tempdima
           \fi
         \fi
         \@tempdimb=\ht1
         \advance\@tempdimb by \dp1
         \advance\@tempdimb by -#2%
         \advance\@tempdimb by #3%
         \leavevmode
         \raise -\@tempdimb \hbox{\box1}%
      \fi
      \egroup%
}%
\def\DFRAME#1#2#3#4#5{%
 \begin{center}
     \let\QCTOptA\empty
     \let\QCTOptB\empty
     \let\QCBOptA\empty
     \let\QCBOptB\empty
     \ifOverFrame 
        #5\QCTOptA\par
     \fi
     \GRAPHIC{#4}{#3}{#1}{#2}{\z@}
     \ifUnderFrame 
        \nobreak\par #5\QCBOptA
     \fi
 \end{center}%
 }%
\def\FFRAME#1#2#3#4#5#6#7{%
 \begin{figure}[#1]%
  \let\QCTOptA\empty
  \let\QCTOptB\empty
  \let\QCBOptA\empty
  \let\QCBOptB\empty
  \ifOverFrame
    #4
    \ifx\QCTOptA\empty
    \else
      \ifx\QCTOptB\empty
        \caption{\QCTOptA}%
      \else
        \caption[\QCTOptB]{\QCTOptA}%
      \fi
    \fi
    \ifUnderFrame\else
      \label{#5}%
    \fi
  \else
    \UnderFrametrue%
  \fi
  \begin{center}\GRAPHIC{#7}{#6}{#2}{#3}{\z@}\end{center}%
  \ifUnderFrame
    #4
    \ifx\QCBOptA\empty
      \caption{}%
    \else
      \ifx\QCBOptB\empty
        \caption{\QCBOptA}%
      \else
        \caption[\QCBOptB]{\QCBOptA}%
      \fi
    \fi
    \label{#5}%
  \fi
  \end{figure}%
 }%
\newcount\dispkind%

\def\makeactives{
  \catcode`\"=\active
  \catcode`\;=\active
  \catcode`\:=\active
  \catcode`\'=\active
  \catcode`\~=\active
}
\bgroup
   \makeactives
   \gdef\activesoff{%
      \def"{\string"}
      \def;{\string;}
      \def:{\string:}
      \def'{\string'}
      \def~{\string~}
    }
\egroup

\def\FRAME#1#2#3#4#5#6#7#8{%
 \bgroup
 \@ifundefined{bbl@deactivate}{}{\activesoff}
 \ifnum\draft=\@ne
   \wasdrafttrue
 \else
   \wasdraftfalse%
 \fi
 \def\LaTeXparams{}%
 \dispkind=\z@
 \def\LaTeXparams{}%
 \doFRAMEparams{#1}%
 \ifnum\dispkind=\z@\IFRAME{#2}{#3}{#4}{#7}{#8}{#5}\else
  \ifnum\dispkind=\@ne\DFRAME{#2}{#3}{#7}{#8}{#5}\else
   \ifnum\dispkind=\tw@
    \edef\@tempa{\noexpand\FFRAME{\LaTeXparams}}%
    \@tempa{#2}{#3}{#5}{#6}{#7}{#8}%
    \fi
   \fi
  \fi
  \ifwasdraft\draft=1\else\draft=0\fi{}%
  \egroup
 }%
\def\TEXUX#1{"texux"}

\def\limfunc#1{\mathop{\rm #1}}%

\long\def\QQQ#1#2{%
     \long\expandafter\def\csname#1\endcsname{#2}}%
\@ifundefined{QTP}{\def\QTP#1{}}{}
\@ifundefined{QEXCLUDE}{\def\QEXCLUDE#1{}}{}
\@ifundefined{Qlb}{}{}
\@ifundefined{Qlt}{}{}
\long\def\QQA#1#2{}%
\def\QTR#1#2{{\csname#1\endcsname #2}}
\def\EXPAND#1[#2]#3{}%
\def\NOEXPAND#1[#2]#3{}%
\def\LaTeXparent#1{}%
\def\ChildStyles#1{}%
\def\ChildDefaults#1{}%
\def\QTagDef#1#2#3{}%
\@ifundefined{StyleEditBeginDoc}{}{}
\def\QQfnmark#1{\footnotemark}

\@ifundefined{INDEX}{\def\INDEX#1#2{}{}}{}%
\@ifundefined{SUBINDEX}{\def\SUBINDEX#1#2#3{}{}{}}{}%
\@ifundefined{initial}%
   {\def\initial#1{\bigbreak{\raggedright\large\bf #1}\kern 2\p@\penalty3000}}%
   {}%
\@ifundefined{entry}{}{}%
\@ifundefined{primary}{}{}%
\@ifundefined{secondary}{}{}%
\@ifundefined{ZZZ}{}{\makeatletter\input gnuindex.sty\makeatother\makeindex\makeatletter}%
\@ifundefined{abstract}{%
 \def\abstract{%
  \if@twocolumn
   \section*{Abstract (Not appropriate in this style!)}%
   \else \small 
   \begin{center}{\bf Abstract\vspace{-.5em}\vspace{\z@}}\end{center}%
   \quotation 
   \fi
  }%
 }{%
 }%
\@ifundefined{endabstract}{\def\endabstract
  {\if@twocolumn\else\endquotation\fi}}{}%
\@ifundefined{maketitle}{\def\maketitle#1{}}{}%
\@ifundefined{affiliation}{\def\affiliation#1{}}{}%
\@ifundefined{proof}{}{}%
\@ifundefined{endproof}{}{}%
\@ifundefined{newfield}{\def\newfield#1#2{}}{}%
\@ifundefined{chapter}{\def\chapter#1{\par(Chapter head:)#1\par }%
 \newcount\c@chapter}{}%
\@ifundefined{part}{\def\part#1{\par(Part head:)#1\par }}{}%
\@ifundefined{section}{\def\section#1{\par(Section head:)#1\par }}{}%
\@ifundefined{subsection}{\def\subsection#1%
 {\par(Subsection head:)#1\par }}{}%
\@ifundefined{subsubsection}{\def\subsubsection#1%
 {\par(Subsubsection head:)#1\par }}{}%
\@ifundefined{paragraph}{\def\paragraph#1%
 {\par(Subsubsubsection head:)#1\par }}{}%
\@ifundefined{subparagraph}{\def\subparagraph#1%
 {\par(Subsubsubsubsection head:)#1\par }}{}%
\@ifundefined{therefore}{}{}%
\@ifundefined{backepsilon}{}{}%
\@ifundefined{yen}{}{}%
\@ifundefined{registered}{%
   \def\registered{\relax\ifmmode{}\r@gistered
                    \else$\m@th\r@gistered$\fi}%
 \def\r@gistered{^{\ooalign
  {\hfil\raise.07ex\hbox{$\scriptstyle\rm\text{R}$}\hfil\crcr
  \mathhexbox20D}}}}{}%
\@ifundefined{Eth}{}{}%
\@ifundefined{eth}{}{}%
\@ifundefined{Thorn}{}{}%
\@ifundefined{thorn}{}{}%
\@ifundefined{degree}{}{}%
\newdimen\theight
\def\Column{%
 \vadjust{\setbox\z@=\hbox{\scriptsize\quad\quad tcol}%
  \theight=\ht\z@\advance\theight by \dp\z@\advance\theight by \lineskip
  \kern -\theight \vbox to \theight{%
   \rightline{\rlap{\box\z@}}%
   \vss
   }%
  }%
 }%
\def\qed{%
 \ifhmode\unskip\nobreak\fi\ifmmode\ifinner\else\hskip5\p@\fi\fi
 \hbox{\hskip5\p@\vrule width4\p@ height6\p@ depth1.5\p@\hskip\p@}%
 }%
\def\miss{\hbox{\vrule height2\p@ width 2\p@ depth\z@}}%
\def\tcol#1{{\baselineskip=6\p@ \vcenter{#1}} \Column}  %
\def\newfmtname{LaTeX2e}
\def\chkcompat{%
   \if@compatibility
   \else
     \usepackage{latexsym}
   \fi
}

\ifx\fmtname\newfmtname
  \DeclareOldFontCommand{\rm}{\normalfont\rmfamily}{\mathrm}
  \DeclareOldFontCommand{\sf}{\normalfont\sffamily}{\mathsf}
  \DeclareOldFontCommand{\tt}{\normalfont\ttfamily}{\mathtt}
  \DeclareOldFontCommand{\bf}{\normalfont\bfseries}{\mathbf}
  \DeclareOldFontCommand{\it}{\normalfont\itshape}{\mathit}
  \DeclareOldFontCommand{\sl}{\normalfont\slshape}{\@nomath\sl}
  \DeclareOldFontCommand{\sc}{\normalfont\scshape}{\@nomath\sc}
  \chkcompat
\fi

\def\alpha{{\Greekmath 010B}}%
\def\beta{{\Greekmath 010C}}%
\def\gamma{{\Greekmath 010D}}%
\def\delta{{\Greekmath 010E}}%
\def\epsilon{{\Greekmath 010F}}%
\def\zeta{{\Greekmath 0110}}%
\def\eta{{\Greekmath 0111}}%
\def\theta{{\Greekmath 0112}}%
\def\iota{{\Greekmath 0113}}%
\def\kappa{{\Greekmath 0114}}%
\def\lambda{{\Greekmath 0115}}%
\def\mu{{\Greekmath 0116}}%
\def\nu{{\Greekmath 0117}}%
\def\xi{{\Greekmath 0118}}%
\def\pi{{\Greekmath 0119}}%
\def\rho{{\Greekmath 011A}}%
\def\sigma{{\Greekmath 011B}}%
\def\tau{{\Greekmath 011C}}%
\def\upsilon{{\Greekmath 011D}}%
\def\phi{{\Greekmath 011E}}%
\def\chi{{\Greekmath 011F}}%
\def\psi{{\Greekmath 0120}}%
\def\omega{{\Greekmath 0121}}%
\def\varepsilon{{\Greekmath 0122}}%
\def\vartheta{{\Greekmath 0123}}%
\def\varpi{{\Greekmath 0124}}%
\def\varrho{{\Greekmath 0125}}%
\def\varsigma{{\Greekmath 0126}}%
\def\varphi{{\Greekmath 0127}}%

\def\nabla{{\Greekmath 0272}}
\def\FindBoldGroup{%
   {\setbox0=\hbox{$\mathbf{x\global\edef\theboldgroup{\the\mathgroup}}$}}%
}

\def\Greekmath#1#2#3#4{%
    \if@compatibility
        \ifnum\mathgroup=\symbold
           \mathchoice{\mbox{\boldmath$\displaystyle\mathchar"#1#2#3#4$}}%
                      {\mbox{\boldmath$\textstyle\mathchar"#1#2#3#4$}}%
                      {\mbox{\boldmath$\scriptstyle\mathchar"#1#2#3#4$}}%
                      {\mbox{\boldmath$\scriptscriptstyle\mathchar"#1#2#3#4$}}%
        \else
           \mathchar"#1#2#3#4%
        \fi 
    \else 
        \FindBoldGroup
        \ifnum\mathgroup=\theboldgroup 
           \mathchoice{\mbox{\boldmath$\displaystyle\mathchar"#1#2#3#4$}}%
                      {\mbox{\boldmath$\textstyle\mathchar"#1#2#3#4$}}%
                      {\mbox{\boldmath$\scriptstyle\mathchar"#1#2#3#4$}}%
                      {\mbox{\boldmath$\scriptscriptstyle\mathchar"#1#2#3#4$}}%
        \else
           \mathchar"#1#2#3#4%
        \fi                 
          \fi}

\newif\ifGreekBold  \GreekBoldfalse
\let\SAVEPBF=\pbf
\def\pbf{\GreekBoldtrue\SAVEPBF}%

\@ifundefined{theorem}{\newtheorem{theorem}{Theorem}}{}
\@ifundefined{lemma}{\newtheorem{lemma}[theorem]{Lemma}}{}
\@ifundefined{corollary}{}{}
\@ifundefined{proposition}{}{}
\@ifundefined{axiom}{}{}
\@ifundefined{remark}{}{}
\@ifundefined{example}{}{}
\@ifundefined{exercise}{}{}
\@ifundefined{definition}{}{}

\@ifundefined{mathletters}{%
  \newcounter{equationnumber}  
  \def\mathletters{%
     \addtocounter{equation}{1}
     \edef\@currentlabel{\theequation}%
     \setcounter{equationnumber}{\c@equation}
     \setcounter{equation}{0}%
     \edef\theequation{\@currentlabel\noexpand\alph{equation}}%
  }
  
}{}

\@ifundefined{BibTeX}{%
    \def\BibTeX{{\rm B\kern-.05em{\sc i\kern-.025em b}\kern-.08em
                 T\kern-.1667em\lower.7ex\hbox{E}\kern-.125emX}}}{}%
\@ifundefined{AmS}%
    {\def\AmS{{\protect\usefont{OMS}{cmsy}{m}{n}%
                A\kern-.1667em\lower.5ex\hbox{M}\kern-.125emS}}}{}%
\@ifundefined{AmSTeX}{}{}%
\begin{document}

\title{Absolutely Rigid Systems and
Absolutely Indecomposable Groups} 
\author{Paul C. Eklof\thanks{Partially supported by NSF Grants
DMS-9501415 and DMS-9704477}\ \ and\ \   
Saharon Shelah\thanks{Partially supported by NSF Grant
DMS-9704477. Pub. No. 678} }
\maketitle

\begin{abstract}
We give a new proof that there are arbitrarily large indecomposable abelian
groups; moreover, the groups constructed are absolutely indecomposable, that
is, they remain indecomposable in any generic extension. However, any
absolutely rigid family of groups has cardinality less than the partition
cardinal $\kappa (\omega )$.
\end{abstract}

\section{Introduction}

Mark Nadel \cite{nadel} asked whether there is a proper class of
torsion-free abelian groups $\{A_{\nu }:\nu \in Ord\}$ with the property
that for any $\nu \neq \mu $, $A_{\nu }$ and $A_{\mu }$ are not $L_{\infty
\omega }$-equivalent; this is the same as requiring that $A_{\nu }$ and $%
A_{\mu }$ do not become isomorphic in any generic extension of the universe.
In that case we say that $A_{\nu }$ and $A_{\mu }$ are absolutely
non-isomorphic \index{absolutely non-isomorphic}. This is not hard to achieve for torsion abelian groups,
since groups of different $p$-length are absolutely non-isomorphic. (See
section 1 for more information.)

Nadel's approach to the question in \cite{nadel} involved looking at known
constructions of rigid systems $\{A_{i}:i\in I\}$ to see if they had the
property that for $i\neq j$, $\limfunc{Hom}(A_{i},A_{j})$ remains zero in
any generic extension of the universe. We call these \textit{absolutely
rigid systems}. \index{absolutely rigid} Similarly we call a group \textit{absolutely rigid} (resp. 
\textit{absolutely indecomposable}) if it is rigid (resp. indecomposable) in
any generic extension. \index{generic extension} Nadel showed that  the Fuchs-Corner construction in 
\cite[\S 89]{F} constructs an absolutely rigid system $\{A_{\nu }:\nu
<2^{\lambda }\}$ of groups of cardinality $\lambda $, where $\lambda $ is
less than the first strongly inaccessible cardinal. But he pointed out that
other constructions, such as Fuchs' construction \cite{F2} of a rigid system
of groups of cardinality the first measurable or Shelah's \cite{Sh74} for an
arbitrary cardinal involve non-absolute notions like direct products or
stationary sets; so the rigid systems constructed may not remain rigid when
the universe of sets is expanded. The same comment applies to any
construction based on a version of the Black Box.

Here we show that there do \textit{not} exist arbitrarily large absolutely
rigid systems. The cardinal $\kappa (\omega )$ in the following theorem is
defined in section 2; it is an inaccessible cardinal much larger than the
first inaccessible, but small enough to be consistent with the Axiom of
Constructibility.

\begin{theorem}
\label{main1} If $\kappa $ is a cardinal $\geq \kappa (\omega )$ and $%
\{A_{\nu }:\nu <\kappa \}$ is a family of non-zero abelian groups, then
there are $\mu \neq \nu $ in $\kappa $ such that in some generic extension $%
V[G]$ of the universe, $V$, there is a non-zero (even one-one) homomorphism $%
f:A_{\nu }\rightarrow A_{\mu }$.
\end{theorem}

This cardinal $\kappa (\omega )$ (called the ``first beautiful cardinal'' by
the second author in \cite{Sh110}) is the precise dividing line:

\begin{theorem}
\label{main2} If $\kappa $ is a cardinal $<\kappa (\omega )$ and $\lambda $
is any cardinal $\geq \kappa (\omega )$, there is a family $\{A_{\mu }:\mu
<\kappa \}$ of torsion-free groups of cardinality $\lambda $ such that in
any generic extension $V[G]$, for all $\mu \in \kappa $, $A_{\mu }$ is
indecomposable and for $\nu \neq \mu $, $\limfunc{Hom}(A_{\nu },A_{\mu })=0$.
\end{theorem}

Despite the limitation imposed by Theorem \ref{main1}, the construction used
to prove Theorem \ref{main2} yields the existence of a proper class of
absolutely different torsion-free groups, in the following strong form. This
answers the question of Nadel in the affirmative, and also provides a new
proof of the existence of arbitrarily large indecomposables.

\begin{theorem}
\label{main3}For each uncountable cardinal $\lambda $, there exist $%
2^{\lambda }$ torsion-free absolutely indecomposable groups $\{H_{i,\lambda
}:i<2^{\lambda }\}$ of cardinality $\lambda $ such that whenever $\lambda
\neq \rho $ or $i\neq j$, $H_{i,\lambda }$ and $H_{j,\rho }$ are absolutely
non-isomorphic.
\end{theorem}

We show that the groups $A_{\mu }$ in Theorem \ref{main2} and the groups $%
H_{i,\lambda }$ in Theorem \ref{main3} are absolutely
 indecomposable by showing that in any generic extension
the only automorphisms they have are $1$ and $-1$. (The proof of Theorem \ref
{main3} does not depend on results from \cite{Sh110}.) However, we cannot
make the groups absolutely rigid:

\begin{theorem}
\label{main4} If $\kappa $ is a cardinal $\geq \kappa (\omega )$ and $A$
is a torsion-free abelian group of cardinality $\kappa$, then
 in some generic extension $%
V[G]$ of the universe, there is an endomorphism of $A$ which
is not multiplication by a  rational number.
\end{theorem}

\section{Infinitary logic and generic extensions}

We will confine ourselves to the language of abelian groups. Thus an \textit{%
atomic formula} is one of the form $\sum_{i=0}^{n}c_{i}x_{i}=0$ where the $%
c_{i}$ are integers and the $x_{i}$ are variables.

$L_{\omega \omega }$ consists of the closure of the atomic formula under
negation ($\lnot $), finite conjunctions ($\wedge $) and disjunctions ($\vee 
$), and existential ($\exists x$) and universal ($\forall x$) quantification
(over a single variable --- or equivalently over finitely many variables). $%
L_{\infty \omega }$ consists of the closure of the atomic formula under
negation, \textit{arbitrary} (possibly infinite) conjunction ($\bigwedge $)
and disjunction ($\bigvee $), and under existential and universal
quantification ($\exists x,\forall x$). Rather than give formal definitions
of other model-theoretic concepts, we will illustrate them with examples.
Thus the formula $\varphi (y):$%
\[
\forall x\exists z(2z=x)\wedge (\lnot \exists z(3z=y)) 
\]
is a formula of $L_{\omega \omega }$ with one free variable, $y$, which
``says'' that every element is $2$-divisible, but $y$ is not divisible by $3$%
. More formally, if $A$ is an abelian group and $a\in A$ we write $A\models
\varphi [a]$ and say ``$a$ satisfies $\varphi $ in $A$'', if and only if
every element of $A$ is divisible by $2$ and there is no $b\in A$ such that $%
3b=a$.

Also, the formula $\psi (x):$%
\[
\exists y(py=x)\wedge (\bigwedge_{n\geq 1}\exists z(p^{n}z=y))\wedge (x\neq
0) 
\]
is a formula of $L_{\infty \omega }$ with free variable $x$ such that $%
A\models \psi [a]$ if and only if $a\in p^{\omega +1}A-\{0\}${.}

A \textit{sentence} is a formula which has no free variables; if $\varphi $
is a sentence of $L_{\infty \omega },$we write $A\models \varphi $ if and
only if $\varphi $ is true in $A$. We write $A\equiv _{\infty \omega }B$ to
mean that every sentence of $L_{\infty \omega }$ true in $A$ is true in $B$
(and conversely because $\lnot \varphi $ true in $A$ implies $\lnot \varphi $
true in $B$.) Obviously, if there is an isomorphism $f:A\rightarrow B$, then 
$A\equiv _{\infty \omega }B$. A necessary and sufficient condition for $%
A\equiv _{\infty \omega }B$ is given by the following (\cite{karp}, or see 
\cite[pp. 13f]{bar}):

\begin{lemma}
\label{b&f}$A\equiv _{\infty \omega }B$ if and only if there is a set $P$ of
bijections $p:A_{p}\rightarrow B_{p}$ from a finite subset $A_{p}$ of $A$
onto a finite subset $B_{p}$ of $B$ with the following properties:

(i)\emph{[the elements of }$P$\emph{\ are partial isomorphisms]} for every
atomic formula $\varphi (x_{1},...,x_{m})$ and elements $a_{1},...,a_{m}$ of 
$\limfunc{dom}(p)$, $A\models \varphi [a_{1},...,a_{m}]$ if and only if $%
B\models \varphi [p(b_{1}),...,p(b_{m})]$;

(ii)\emph{[the back-and-forth property]} for every $p\in P$ and every $a\in A
$ (resp. $b\in B$), there is $p^{\prime }\in P$ such that $p\subseteq
p^{\prime }$ and $a\in \limfunc{dom}(p^{\prime })$ (resp. $b\in \limfunc{rge}%
(p^{\prime })$).
\end{lemma}

It is an easy consequence that if $A$ and $B$ are countable, then $A\equiv
_{\infty \omega }B$ if and only if $A\cong B$. Also, this implies that if it
is true in $V$ that $A\equiv _{\infty \omega }B$, then $A\equiv _{\infty
\omega }B$ remains true in $V[G]$. The converse is easy, and direct, since a
sentence of $L_{\infty \omega }$ which, in $V$, holds true in $A$ but false
in $B$ has the same status in $V[G]$, since no new elements are added to the
groups.

(By a generic extension we mean an extension of the universe $V$ of sets
defined by the method of forcing. In general,  more sets are added to the
universe; possibly, for example,  a bijection between an uncountable
cardinal $\lambda $  and the countable set $\omega $. So cardinals of $V$ 
may not be cardinals in $V[G]$; but the ordinals of $V[G]$ are the same as
the ordinals of $V$. Also, the elements of any set in $V$ are the same in $V$
or $V[G]$.)

There exist non-isomorphic uncountable groups $A$ and $B$ (of cardinality $%
\aleph _{1}$ for example) such that $A\equiv _{\infty \omega }B$. (See for
example \cite{Ek1}.) However, for any groups $A$ and $B$ in the universe, $V$%
, there is a generic extension $V[G]$ of $V$ in which $A$ and $B$ are both
countable (cf. \cite[Lemma 19.9, p. 182]{J}). Therefore we can conclude that 
$A\equiv _{\infty \omega }B$ if and only if $A$ and $B$ are ``potentially
isomorphic'', that is, there is a generic extension $V[G]$ of the universe
in which they become isomorphic. Barwise argues in \cite[p. 32]{bar} that
potential isomorphism (that is, the relation $\equiv _{\infty \omega }$) is
``a very natural notion of isomorphism, one of which mathematicians should
be aware. If one proves that $A\ncong B$ but leaves open the question [of
whether $A$ and $B$ are potentially isomorphic] then one leaves the
possibility that $A$ and $B$ are not isomorphic for trivial reasons of
cardinality. Or to put it the other way round, a proof that [$A$ is not
potentially isomorphic to $B$ ] is a proof that $A\ncong B$ for nontrivial
reasons.''

As an example, consider reduced $p$-groups $A_{\nu }$ ($\nu $ any ordinal)
such that the length of $A_{\nu }$ is $\nu $, that is, $p^{\nu }A_{\nu }=0$
but for all $\mu <\nu $, $p^{\mu }A_{\nu }\neq 0$. Then for any ${\nu }%
_{1}\neq {\nu }_{2}${,} the groups $A_{{\nu }_{1}}$ and $A_{{\nu }_{2}}$ are
not even potentially isomorphic: this is because for any $\nu $ there is a
formula $\theta _{\nu }(x)$ such that $\exists x(\theta _{\nu }(x)\wedge
x\neq 0)$ is true in a $p$-group $A$ if and only if $A$ has length $\geq \nu 
$. Indeed, we define, by induction, $\theta _{\nu }$ to be 
\[
\exists y(py=x\wedge \theta _{\mu }(y)) 
\]
if $\nu =\mu +1$ and if $\nu $ is a limit ordinal define $\theta _{\nu }$ to
be 
\[
\bigwedge_{\mu <\nu }\exists y(py=x\wedge \theta _{\mu }(y))\text{.} 
\]

Thus there is a proper class of pairwise absolutely non-isomorphic $p$%
-groups. Although there is not available a standard group-theoretic notion
which will serve the same purpose for torsion-free groups, we will prove in
section 5 that there is a proper class of indecomposable torsion-free
abelian groups $\{H_{{\lambda }}:{\lambda }$ a cardinal$\}$ such that for
any ${\lambda }\neq \rho ${,} the groups $H_{\lambda }$ and $H_{\rho }$ are
not $L_{\infty \omega }$-equivalent.

\section{Quasi-well-orderings and beautiful cardinals}

A \textit{quasi-order} $Q$ is a pair $(Q,\leq _{Q})$ where $\leq _{Q}$ is a
reflexive and transitive binary relation on $Q$. There is an extensive
theory of well-orderings of quasi-orders developed by Higman, Kruskal,
Nash-Williams and Laver among others (cf. \cite{nw}, \cite{laver}). A
generalization to uncountable cardinals is due to the second author (\cite
{Sh110}). The key notion that we need is the following: for an infinite
cardinal $\kappa $, $Q$ is called $\kappa $\textit{--narrow } if there is no 
\textit{antichain} in $Q$ of size $\kappa $, i.e., for every $f:\kappa
\rightarrow Q$ there exist $\nu \neq \mu $ such that $f(\nu )\leq _{Q}f(\mu )
$. (Note that this use of the terminology ``antichain'' --- in \cite[p.32]
{lav2} for example --- is different from its use in forcing theory.)

A\textit{\ tree} \index{tree} is a partially-ordered set $(T,\leq )$ such that for all $%
t\in T$, pred($t$)$=\{s\in T:s<t\}$ is a well-ordered set; moreover, there
is only one element $r$ of $T$, called the \textit{root} of $T$, such that
pred($r$) is empty. The order-type of pred($t$) is called the height of $t$,
denoted ht$(t)$; the height of $T$ is $\sup \{$ht$(t) + 1:t\in T\}$.

If $Q$ is a quasi-order, \index{quasi-order} a $Q$\textit{-labeled tree} is a pair $(T,\Phi
_{T}) $ consisting of a tree $T$ of height $\leq \omega $ and a function $%
\Phi _{T}:T\rightarrow Q$. On any set of $Q$-labeled trees we define a
quasi-order by: $(T_{1},\Phi _{1})\preceq (T_{2},\Phi _{2})$ if and only if
there is a function $\theta :T_{1}\rightarrow T_{2}$ which preserves the
tree-order (i.e. $t\leq _{T_{1}}t^{\prime }$ implies $\theta (t)\leq
_{T_{2}}\theta (t^{\prime })$) as well as the height of elements and also is
such that for all $t\in T_{1}$, $\Phi _{1}(t)\leq _{Q}\Phi _{2}(\theta (t))$.

One result from \cite{Sh110} that we will use implies that for sufficiently
large cardinals $\kappa $ and sufficiently small $Q$, any set of $Q$-labeled
trees is $\kappa $-narrow. In order to state the result precisely we need to
define a certain (relatively small) large cardinal.

Let $\kappa (\omega )$ be the first $\omega $-Erd\"{o}s cardinal, 
\index{Erd\"{o}s cardinal} i.e., the
least cardinal such that $\kappa \longrightarrow (\omega )_{{}}^{<\omega }$;
in other words, the least cardinal such that for every function $F$ from the
finite subsets of $\kappa $ to $2$ there is an infinite subset $X$ of $%
\kappa $ such that there is a function $c:\omega \rightarrow 2$ such for
every finite subset $Y$ of $X$, $F(Y)=c(|Y|)$. It has been shown that this
cardinal is strongly inaccessible (cf. \cite[p. 392]{J}). Thus it cannot be
proved in ZFC that $\kappa (\omega )$ exists (or even that its existence is
consistent). If it exists, there are many weakly compact cardinals below it,
and, on the other hand, it is less than the first measurable cardinal (if
such exists). Moreover, if it is consistent with ZFC that there is such a
cardinal, then it is consistent with ZFC + V = L that there is such a
cardinal (\cite{sil}). If $\kappa (\omega )$ does not exist, then Theorem \ref{qwo1} is
uninteresting. On the other hand, Theorem \ref{qwo2} then applies to every
cardinal $\kappa $, and its consequences, given in section 4, are still of
interest.

The following is a consequence of results proved in \cite{Sh110} (cf.
Theorem 5.3, p. 208 and Theorem 2.10, p. 197):

\begin{theorem}
\label{qwo1}If $Q$ is a quasi-order of cardinality $<\kappa (\omega )$, and $%
\mathcal{S}$ is a set of $Q$-labeled trees, then $\mathcal{S}$ is $\kappa
(\omega )$-narrow.
\end{theorem}

On the other hand, it follows from results in \cite{Sh110} that for any
cardinal smaller than $\kappa (\omega )$, there is an absolute antichain of
that size:

\begin{theorem}
\label{qwo2}If $\kappa <\kappa (\omega )$, there is a family $\mathcal{T}%
=\{(T_{{\mu }},\Phi _{{\mu }}):{\mu }<\kappa \}$ of $\omega $-labeled trees,
each of cardinality $<\kappa (\omega )$, such that in any generic extension
of $V$, for all $\mu \neq \nu $, $(T_{{\mu }},\Phi _{{\mu }})\npreceq
(T_{\nu },\Phi _{\nu })$.
\end{theorem}

Some commentary is needed on the absoluteness of the antichain $\mathcal{T}%
=\{(T_{{\mu }},\Phi _{{\mu }}):{\mu }<\kappa \}$, since this is not dealt
with directly in \cite{Sh110}. $\mathcal{T}$ is constructed in a concrete,
absolute, way from a function $F$ which is an example witnessing the fact
that $\kappa <\kappa (\omega )$. First, a $\kappa $-$D$-barrier $B$ and
function $q:B\rightarrow \omega $ is constructed ($B$ is a kind of elaborate
indexing for an antichain cf. \cite[proof of 2.5, p. 195]{Sh110}). This
gives rise (\cite[proof of 1.12, pp. 192f]{Sh110}) to an example showing
that $\mathcal{P}_{\beta }(\omega )$ is not $\kappa $-narrow for some $\beta
<\kappa (\omega )$; this example is embedded into the quasi-order of $\omega 
$-labeled trees, giving rise to $\mathcal{T}$ (\cite[p. 221]{Sh110}). The
proof that $\mathcal{T}$ is an antichain reduces to the key property of $F$,
a property which is absolute by an argument of Silver \cite{sil}. Using an
equivalent definition of $\kappa (\omega )$, $F$ is taken to be a function
from the finite subsets of $\kappa $ to $\omega $ such that there is no
one-one function $\sigma :\omega \rightarrow \omega $ such that for all $%
n\in \omega $, $F(\{\sigma (0),...,\sigma (n-1)\})=F(\{\sigma (1),...,\sigma
(n)\})$; this property of $F$ is preserved under generic extensions because
it is equivalent to the well-foundedness of a certain tree; more precisely,
the tree of finite partial attempts at $\sigma $ has no infinite branch.

\section{A bound on the size of absolutely rigid systems}

In this section we will prove Theorems \ref{main1} and \ref{main4}.. Suppose
that $\{A_\nu :\nu <\kappa \}$ is a family of non-zero abelian groups, where
we can assume that $\kappa =\kappa (\omega )$. For each $\nu <\kappa $, let $%
T_\nu $ be the tree of finite sequences of elements of $A_\nu $; that is,
the elements of $T_\nu $ are 1-1 functions $s:n_s\rightarrow A_\nu $ for
some $n_s\in \omega $ and $s\leq t$ if and only if $n_s\leq n_t$ and $%
t\upharpoonright n_s=s$.

Let $Q_{ab}$ be the set of all quantifier-free $n$-types of abelian groups;
that is, $Y\in Q_{ab}$ if and only if for some abelian group $G$, some $n\in
\omega $, and some function $s:n\rightarrow G$, $Y$ is the set $\limfunc{tp}%
_{qf}(s/G)$ of all quantifier-free formulas $\varphi (x_{0},...,x_{n-1})$ of 
$L_{\omega \omega }$ such that $G\vDash \varphi [s(0),...,s(n-1)]$.
Partially-order $Q_{ab}$ by the relation of inclusion.

Define $\Phi _{\nu }:T_{\nu }\rightarrow Q_{ab}$ by letting $\Phi _{\nu }(s)=%
\limfunc{tp}_{qf}(s/A_{\nu })$. Now we can apply Theorem \ref{qwo1} to the
family of $Q_{ab}$-labeled trees $\mathcal{S}=\{(T_{\nu },\Phi _{\nu }):\nu
<\kappa \}$. (Note that the cardinality of $Q_{ab}$ is $2^{\aleph _{0}}$
which is $<\kappa (\omega )$ since $\kappa (\omega )$ is
strongly-inaccessible.) Therefore there exists $\nu \neq \mu $ such that $%
(T_{\nu },\Phi _{\nu })\preceq (T_{\mu },\Phi _{\mu })$, say $\theta :T_{\nu
}\rightarrow T_{\mu }$ is such that $s\leq t$ implies $\theta (s)\leq \theta
(t)$ and for all $s\in T_{\nu }$, $\Phi _{\nu }(s)\subseteq \Phi _{\mu
}(\theta (s))$.

Now move to a generic extension $V[G]$ in which $A_{\nu }$ is countable. In $%
V[G]$, let $\sigma :\omega \rightarrow A_{\nu }$ be a surjection. We will
define an embedding $f:A_{\nu }\rightarrow A_{\mu }$ by letting $f(\sigma
(n))=\theta (\sigma \upharpoonright n+1)(n)$ for all $n<\omega $. To see
that $f$ is an embedding, note that $f(\sigma (n))=\theta (\sigma
\upharpoonright k)(n)$ for all $k>n$ since $\theta $ preserves the tree
ordering; moreover, for any $a,b,c\in A_{\nu }$, there is a $k$ such that $%
a,b,c\in \limfunc{rge}(\sigma \upharpoonright k)$ so since 
$$
\Phi _{\nu }(\sigma \upharpoonright k)\subseteq \Phi _{\mu }(\theta (\sigma
\upharpoonright k))=\Phi _{\mu }(\left\langle f(0),...,f(k-1)\right\rangle ) 
$$
every quantifier-free formula satisfied by $a,b,c$ in $A_{\nu }$ (e.g. $%
a\neq 0$, $a-b=c$, $ab=c$) is satisfied by $f(a),f(b),f(c)$ in $A_{\mu }$.
This completes the proof of Theorem \ref{main1}

\medskip

The argument is very general and could be applied to any family of
structures (for example, to those in \cite{GM}). If we start with a
torsion-free group $A$ of cardinality $\kappa \geq \kappa (\omega )$, and
apply the argument to the family of structures $\{\left\langle
A,a_v\right\rangle :\nu <\kappa (\omega )\}$ where $\{a_\nu :\nu <\kappa
(\omega )\}$ is a linearly independent subset of $A$, then we obtain $\nu
\neq \mu $ such that in a generic extension in which $A$ becomes countable
we have an embedding  $f:A\rightarrow A$ taking $a_\nu $ to $a_\mu $. This
proves Theorem \ref{main4}.

Ernest Schimmerling has pointed out that there is a ``soft'' proof of these
results (not relying on Theorem \ref{qwo1})  using a model of set theory
(with $\{A_\nu :\nu <\kappa \}$ as additional predicate) and a set of
indiscernibles given by the defining property of $\kappa =\kappa (\omega )$.

\section{Existence theorem}

In this section we will prove Theorem \ref{main2}. So let $\kappa <\kappa
(\omega )$ and let $\lambda $ be a cardinal $\geq \kappa (\omega )$. Let $%
\{(T_{{\mu }},\Phi _{{\mu }}):{\mu }<\kappa \}$ be the family of $\omega $%
-labeled trees as in Theorem \ref{qwo2}. We can assume that each node of $T_{%
{\mu }}$ of height $m$ is a sequence of length $m$ and the tree-ordering is
extension of sequences; so the root of the tree is the empty sequence $<>$.

Let $\left\langle p_{n,m,j}:n,m\in \omega ,j\in \{0,1\}\right\rangle $ and $%
\left\langle q_{n,m,\ell ,j}:n,m,\ell \in \omega ,j\in \{0,1\}\right\rangle $
be two lists, with no overlap, of distinct primes.

For any ordinal $\alpha $, $Z_{\alpha }$ be the tree of finite strictly
decreasing non-empty sequences $z$ of ordinals $\leq \alpha $ such that $%
z(0)=\alpha $. (Thus for some $n\in \omega $, $z:n\rightarrow \alpha $ such
that $\alpha =z(0)>z(1)>...>z(n-1)$.)

For $n\in \omega $ let $g_{n}:\lambda \rightarrow \mathcal{P}([\lambda
n,\lambda (n+1)))$ such that for each $\nu <\lambda $, $g_{n}(\nu )$ is a
subset of $[\lambda n,\lambda (n+1))$ which is cofinal in $\lambda
(n+1)=\lambda n+\lambda $. (Here the operations are ordinal addition and
multiplication, so, in particular, $\lambda n$ is less than $\lambda ^{+}$,
the cardinal successor of $\lambda $.) We also require that for $\mu \neq
\nu $, $g_{n}(\mu )\cap g_{n}(\nu )=\emptyset $. For $n>0$, let $%
Y_{n}=\bigcup \limfunc{rge}(g_{n})$, and let $Y_{0}=g_{0}(0)$.

For ${\mu }<\kappa $, let $W_{{\mu }}$ be the $\Bbb{Q}$-vector space with
basis $\bigcup_{n\in \omega }\mathcal{A}_{n}\cup \mathcal{B}_{n,\mu }$ where
for $n>0$%
\begin{eqnarray*}
\mathcal{A}_{n} &=&\{a_{z}^{\alpha }:\alpha \in Y_{n-1},z\in Z_{\alpha }\}%
\text{, } \\
\mathcal{B}_{n,\mu } &=&\{b_{\eta ,\mu }^{\alpha }:\alpha \in Y_{n-1},\eta
\in T_{{\mu }}-\{<>\}\}
\end{eqnarray*}
and $\mathcal{A}_{0}=\{a_{\mu }^{0}\}=\mathcal{B}_{0,\mu }$. We are going to
define $A_{{\mu }}$ to be a subgroup of $W_{\mu }$. Since $\mu $ is fixed
throughout the construction, we will usually omit the subscript $\mu $ from
what follows (until we come to consider $\limfunc{Hom}(A_{\nu },A_{\mu })$).

For each $n>0$, let $h_{n}$ be a bijection from $\mathcal{A}_{n}\cup 
\mathcal{B}_{n}$ onto $\lambda $; let $h_{0}(a^{0})=0$. Then for any $w\in 
\mathcal{A}_{n}\cup \mathcal{B}_{n}$, and any $\alpha \in g_{n}(h_{n}(w))$,
we will use $a_{<>}^{\alpha }$ or $b_{<>}^{\alpha }$ as a notation for $w$.
(So $a_{<>}^{\alpha }=b_{<>}^{\alpha }$; moreover, $a_{<>}^{\alpha
}=a_{<>}^{\beta }$ if and only if $\alpha $ and $\beta $ belong to the same
member of the range of $g_{n}$.) Now we can define $A$ (= $A_{\mu }$) to be
the subgroup of $W$ generated (as abelian group) by the union of 
\begin{equation}
\bigcup_{n\geq 0}\{\frac{1}{p_{n,m,0}^{k}}{}a_{z}^{\alpha }:{}m,k\in \omega
,z\in Z_{\alpha }\cup \{<>\},\alpha \in Y_{n},\limfunc{dom}(z)=m\}
\label{set1}
\end{equation}

\begin{equation}
\bigcup_{n\geq 0}\{\frac{1}{p_{n,m,1}^{k}}(a_{z}^{\alpha
}+a_{z\upharpoonright m-1}^{\alpha }){:}m,k\in \omega -\{0\},z\in Z_{\alpha
},\alpha \in Y_{n},\limfunc{dom}(z)=m\}  \label{set2}
\end{equation}

and

\begin{equation}
\bigcup_{n\geq 0}\{\frac{1}{q_{n,m,\ell }^{k}}(b_{\eta }^{\alpha }+b_{\eta
\upharpoonright m-1}^{\alpha }){:}m,k\in \omega -\{0\},\eta \in T_{\mu
},\alpha \in Y_{n},\limfunc{dom}(\eta )=m,\Phi _{\mu }(\eta )=\ell \}\text{.}
\label{set3}
\end{equation}

\noindent (where $b_{\eta \upharpoonright -1}^{\alpha }=0$). We will use the sets (\ref{set1}) and (\ref{set2}) to prove that (I) $A_{\mu
}$ is absolutely indecomposable and the set (\ref{set3}) to
prove that (II) $\limfunc{Hom}(A_{\mu },A_{\nu })=0$ for all $\mu \neq \nu $.

If $x\in A_{\mu }$, we will write $p_{{}}^{\infty }|x$ if for every $k\in
\omega $, there exists $v\in A_{\mu }$ such that $p^{k}v=x$. For example, if 
$w\in \mathcal{A}_{n}\cup \mathcal{B}_{n}$ and $\alpha \in g_{n}(h_{n}(w))$
and $\Phi _{\mu }(<>)=\ell _{0}$, then $p_{n,0,0}^{\infty }|w$ and $%
q_{n,0,\ell _{o}}^{\infty }|w$. Assertions about divisibility in $A_{\mu }$
are easily checked by considering the coefficients of linear combinations
over $\Bbb{Q}$ of elements of the basis $\bigcup_{n\in \omega }\mathcal{A}%
_{n}\cup \mathcal{B}_{n}$ of $W_{\mu }$; for example, $p_{n,m,0}^{\infty }|x$
if and only if $x=\sum_{i=1}^{r}c_{i}a_{z_{i}}^{\alpha _{i}}$ for some $%
\alpha _{1},...,\alpha _{r}$ in $Y_{n}$, $z_{i}$ of length $m$, and $%
c_{i}\in \Bbb{Q}$ (with denominator a power of $p_{n,m,0}$).

\smallskip\ 

\begin{center}
\textbf{Proof of (I)}
\end{center}

We will show, in fact, that in any generic extension $V[G]$ the only
automorphisms of $A$ ($=A_{\mu }$) are the trivial ones, $1$ and $-1$. This
part of the proof does not use the trees in $\mathcal{T}$; the absoluteness
is a consequence of an argument using formulas of $L_{\infty \omega }$,
which therefore works in any generic extension. We will use the following
claim:

\begin{quote}
(1A) there are formulas $\psi _{n,\alpha }(x)$ of $L_{\infty \omega }$ ($%
n\in \omega $, $\alpha \in Y_{n}$) such that for any $u\in A$, $A\models
\psi _{n,\alpha }[u]$ if and only if there are $w_{1},...,w_{r}\in \mathcal{A%
}_{n}\cup \mathcal{B}_{n}$, and $c_{1},...,c_{r}\in \Bbb{Z}-\{0\}$, such
that $u=\sum_{i=1}^{r}c_{i}w_{i}$ and $\alpha \in
\bigcup_{i=1}^{r}g_{n}(h_{n}(w_{i}))$.
\end{quote}

Assuming the claim for now, suppose that in $V[G]$ there is an automorphism $%
F$ of $A$. For any $n\in \omega $, consider any $w\in \mathcal{A}_{n}\cup 
\mathcal{B}_{n}$; since $w=a_{<>}^{\alpha }$ for $\alpha \in g_{n}(h_{n}(w))$%
, $p_{n,0,0}^{\infty }|w$; therefore $p_{n,0,0}^{\infty }|F(w)$, and hence $%
F(w)=\sum_{i=1}^{r}c_{i}w_{i}$ for some distinct $w_{i}\in \mathcal{A}%
_{n}\cup \mathcal{B}_{n}$. Moreover, by (1A), $A\models \psi _{n,\alpha }[w]$
if and only if $\alpha \in g_{n}(h_{n}(w))$ if and only if $A\models \psi
_{n,\alpha }[F(w)]$. Thus, since the elements of the range of $g_{n}$ are
disjoint, we must have that $r=1$ and $w_{1}=w$, that is, $F(w)=cw$ for some 
$c=c(w)\in \Bbb{Q}$.

If we can show that $c(w)=c(a^{0})$ for all $w\in \bigcup_{n\in \omega }%
\mathcal{A}_{n}\cup \mathcal{B}_{n}$, then $F$ is multiplication by $%
c(a^{0}) $, and then it is easy to see that $c(a^{0})$ must be $\pm 1$. It
will be enough to show that if $w=a_{z}^{\alpha }$ (resp. $w=b_{\eta
}^{\alpha }$) for some $\alpha \in Y_{n-1}$, then $c(w)=c(a_{<>}^{\alpha })$%
, for then $c(w)=c(w^{\prime })$ for some $w^{\prime }\in \mathcal{A}%
_{n-1}\cup \mathcal{B}_{n-1}$ (namely, the unique $w^{\prime }$ such that $%
\alpha \in g_{n-1}(h_{n-1}(w^{\prime }))$) and by induction $c(w^{\prime
})=c(a^{0})$.

So suppose $w=a_{z}^{\alpha }$; the proof will be by induction on the length
of $z$ that $c(a_{z}^{\alpha })=c(a_{<>}^{\alpha })$. Suppose that the
length of $z=m>0$. Let $c=c(a_{z}^{\alpha })$ and $c^{\prime
}=c(a_{z\upharpoonright m-1}^{\alpha })$. By induction it is enough to prove
that $c=c^{\prime }$. Since $p_{n,m,1}^{\infty }|(a_{z}^{\alpha
}+a_{z\upharpoonright m-1}^{\alpha })$, it is also the case that $%
p_{n,m,1}^{\infty }$ divides 
\[
F(a_{z}^{\alpha })+F(a_{z\upharpoonright m-1}^{\alpha })=ca_{z}^{\alpha
}+c^{\prime }a_{z\upharpoonright m-1}^{\alpha }=c(a_{z}^{\alpha
}+a_{z\upharpoonright m-1}^{\alpha })+(c^{\prime }-c)a_{z\upharpoonright
m-1}^{\alpha } 
\]
so $p_{n,m,1}^{\infty }|(c^{\prime }-c)a_{z\upharpoonright m-1}^{\alpha }$,
which is impossible unless $c=c^{\prime }\,$.

The proof is similar if $w=b_{\eta }^{\alpha }$, but uses the primes $%
q_{n,m,\ell }$. So it remains to prove (1A). We will begin by defining some
auxiliary formulas of $L_{\infty \omega }$. First, we will use $p^{\infty
}|x $ as an abbreviation for 
\[
\bigwedge_{k\in \omega }\exists v_{k}(p^{k}v_{k}=x)\text{.} 
\]
Define $\varphi _{n,m,0}(y)$ to be $p_{n,m,0}^{\infty }|y$. Then for $u\in A$%
, $A\models \varphi _{n,m,0}[u]$ if and only if $u$ is in the subgroup ($\Bbb%
{Z}$-submodule) generated by $\{\frac{1}{p_{n,m,0}^{k}}a_{z}^{\alpha
}:\alpha \in Y_{n},k\in \omega ,z\in Z_{\alpha },\limfunc{dom}(z)=m\}$.
Define $\varphi _{n,m,\beta }(y)$ for each $m>0$ by recursion on $\beta $:
if $\beta =\gamma +1$, $\varphi _{n,m,\beta }(y)$ is 
\[
\varphi _{n,m,\gamma }(y)\wedge \exists y^{\prime }(\varphi _{n,m+1,\gamma
}(y^{\prime })\wedge (p_{n,m+1,1}^{\infty }|(y+y^{\prime }))\text{.} 
\]
If $\beta $ is a limit ordinal, let $\varphi _{n,m,\beta }(y)$ be 
\[
\bigwedge_{\gamma <\beta }\varphi _{n,m,\gamma }(y)\text{.} 
\]

\noindent Then for $u\in A$, $A\models \varphi _{n,m,\beta }[u]$ if and only
if $u$ is in the subgroup generated by 
\[
\{a_{z}^{\alpha }:\alpha \in Y_{n},z\in Z_{\alpha },\limfunc{dom}(z)=m\text{
and }z(m-1)\geq \beta \}\text{.} 
\]
In particular, for $m=1$, recalling that $z\in Z_{\alpha }$ satisfies $%
z(0)=\alpha $, we have that $A\models \varphi _{n,1,\beta }[u]$ if and only
if $u$ is in the subgroup generated by 
\[
\{a_{<\alpha >}^{\alpha }:\alpha \in Y_{n},\alpha \geq \beta \}\text{.} 
\]
Now define $\psi _{n,\alpha }(x)$ to be 
\[
\varphi _{n,0,0}(x)\wedge \exists y[p_{n,1,1}^{\infty }|(x+y)\wedge \varphi
_{n,1,\alpha }(y)\wedge \lnot \varphi _{n,1,\alpha +1}(y)]\text{.} 
\]
If $u=\sum_{i=1}^{r}c_{i}w_{i}$, for some $w_{i}\in \mathcal{A}_{n}\cup 
\mathcal{B}_{n}$, then $p_{n,1,1}^{\infty }|(u+y)$ iff $y=%
\sum_{i=1}^{r}c_{i}a_{<\alpha _{i}>}^{\alpha _{i}}$ for some $\alpha _{i}\in
g_{n}(h_{n}(w_{i}))$; using the cofinality of members of the range of $g_{n}$%
, it follows easily that $\psi _{n,\alpha }(x)$ has the desired property.

\smallskip\ 

\begin{center}
\textbf{Proof of (II)}
\end{center}

Suppose that there is a non-zero homomorphism $H:A_{\nu }\rightarrow A_{\mu }
$. We are going to use $H$ to define $\theta :T_{\nu }\rightarrow T_{\mu }$
showing that $(T_{\nu },\Phi _{\nu })\preceq (T_{\mu },\Phi _{\mu })$,
contrary to the choice of the family of labeled trees. Now $H(w)\neq 0$ for
some $w\in \mathcal{A}_{n}\cup \mathcal{B}_{n,\nu }$ for some $n\in \omega $%
. Thus for some $\alpha \in Y_{n}$, $H(b_{<>,\nu }^{\alpha })\neq 0$. Fix
such an $\alpha $ (which can in fact be any member of $g_{n}(h_{n}(w))$).
Let $\Phi _{\nu }(<>\nolinebreak)=\ell _{0}$. Then $q_{n,0,\ell _{0},0}^{\infty
}|b_{<>,\nu }^{\alpha }$ so $q_{n,0,\ell _{0},0}^{\infty }|H$ $(b_{<>,\nu
}^{\alpha })$, and hence $H(b_{<>,\nu }^{\alpha })$ must be of the form $%
\sum_{i=1}^{r}c_{i}b_{<>,\mu }^{\alpha _{i}^{{}}}$ where $c_{i}\in \Bbb{Q}%
-\{0\}$, and $\Phi _{\mu }(<>)=\ell _{0}$. So letting $\theta (<>)= \ <>$ (as
it must), we have confirmed that $\Phi _{\mu }(\theta (\eta ))=\Phi _{\nu }({%
\eta })$ for ${\eta }=<>$.

Now suppose that for some $m\geq 0$, $\theta ({\eta })$ has been defined for
all nodes ${\eta }$ of $T_{\nu }$ of height $\leq m$ such that $\Phi _{\mu
}(\theta ({\eta }))=\Phi _{\nu }({\eta })$. Moreover, suppose that for every 
$\eta $ of height $\leq m$, the coefficient of $b_{\theta (\eta ),\mu
}^{\alpha }$ in $H(b_{\eta ,\nu }^{\alpha })$ is non-zero. Now consider any
node $\zeta $ of $T_{\nu }$ of height $m+1$; let $\eta =\zeta
\upharpoonright m$. In $A_{\nu }$, for $\ell =\Phi _{\nu }(\zeta )$, $%
q_{n,m+1,\ell ,1}^{\infty }|b_{\eta ,\nu }^{\alpha }+b_{\zeta ,\nu }^{\alpha
}$ so $q_{n,m+1,\ell ,1}^{\infty }|H$ $(b_{\eta ,\nu }^{\alpha })+H$ $%
(b_{\zeta ,\nu }^{\alpha })$ in $A_{\mu }$. Since the coefficient --- call
it $c$ --- of $b_{\theta (\eta ),\mu }^{\alpha }$ in $H(b_{\eta ,\nu
}^{\alpha })$ is non-zero, there must be a node $\zeta ^{\prime }$ in $%
T_{\mu }$ of height $m+1$ such that $\zeta ^{\prime }\upharpoonright
m=\theta (\eta )$ and the coefficient of $b_{\zeta ^{\prime },\mu }^{\alpha
} $ in $H(b_{\zeta ,\nu }^{\alpha })$ is $c$, and moreover such that $\Phi
_{\mu }(\zeta ^{\prime })=\ell $. So we can let $\theta (\zeta )=\zeta
^{\prime }$. This completes the proof of Theorem \ref{main2}.

  \section{Absolutely non-isomorphic indecomposables}

In this section we sketch how to modify the construction in the preceding
section in order to prove Theorem \ref{main3}. (Note that this construction
does not require the trees $T_{\mu }$ of Theorem \ref{qwo2}.) Let $%
\left\langle p_{n,m,j}:n,m\in \omega ,j\in \{0,1\}\right\rangle $ be a list
of distinct primes. Fix an uncountable $\lambda $; for any $\alpha <\lambda
\omega $, let $Z_{\alpha }$ be defined as before. Let $\left\langle
S_{i,\lambda }:i<2^{\lambda }\right\rangle $ be a list of $2^{\lambda }$
distinct subsets of $\lambda $, each of cardinality $\lambda $ (and hence
cofinal in $\lambda $).

For $n\neq 1$, let $g_{n}:\lambda \rightarrow \mathcal{P}([\lambda n,\lambda
(n+1)))$ be defined as before. For $i<2^{\lambda }$, define $g_{1,i}:\lambda
\rightarrow \mathcal{P}([\lambda ,\lambda +\lambda )))$ as before but with
the additional stipulation that for all $\nu <\lambda $, $g_{1,i}(\nu
)\subseteq \{\lambda +\gamma :\gamma \in S_{i,\lambda }\}$. (Here again the
operation is ordinal addition.) Let $Y_{1,i}=\bigcup \limfunc{rge}(g_{1,i})$
we will also choose $g_{1,i}$ such that $Y_{1,i}=\{\lambda +\gamma :\gamma
\in S_{i,\lambda }\}$. For convenience, for $n\neq 1$ we let $Y_{n,i}$
denote $Y_{n}$ (independent of $i$).

For each $n>0$, let $h_{n,i}$ be a bijection from $\{a_{z}^{\alpha }:\alpha
\in Y_{n-1,i}$, $z\in Z_{\alpha }\}$ onto $\lambda $; use these bijections
to make identifications as in the previous construction.

Then $H_{i,\lambda }\ $ is defined to be the subgroup of the $\Bbb{Q}$%
-vector space with basis 
\[
\{a^{0}\}\cup \{a_{z}^{\alpha }:n>0\text{, }\alpha \in Y_{n-1,i}\text{, }%
z\in Z_{\alpha }\} 
\]
which is generated (as abelian group) by the union of 
\[
\bigcup_{n\geq 0}\{\frac{1}{p_{n,m,0}^{k}}{}a_{z}^{\alpha }:{}m,k\in \omega
,z\in Z_{\alpha }\cup \{<>\},\alpha \in Y_{n,i},\limfunc{dom}(z)=m\} 
\]
and

\[
\bigcup_{n\geq 0}\{\frac 1{p_{n,m,1}^k}(a_z^\alpha +a_{z\upharpoonright
m-1}^\alpha ){:}m,k\in \omega -\{0\},z\in Z_\alpha ,\alpha \in Y_{n,i},%
\limfunc{dom}(z)=m\}\text{.} 
\]

As before, the groups $H_{i,\lambda }$ are absolutely indecomposable. It
remains to show that for $\lambda \neq \rho $ or $i\neq j$, $H_{i,\lambda }$
and $H_{j,\rho }$ are not $L_{\infty \omega }$-equivalent (and hence not
isomorphic in any generic extension). For this we use the formulas $\psi
_{1,\alpha }(x)$. If $\lambda =\rho $ and $i\neq j$, without loss of
generality there exists $\gamma \in S_{i,\lambda }-S_{j,\lambda }$; let $%
\alpha =\lambda +\gamma $. If $\lambda <\rho $, let $\alpha =\lambda +\gamma 
$ for any $\gamma $ in any $S_{i,\lambda }$. In either case, $\exists x\psi
_{1,\alpha }(x)$ is true in $H_{i,\lambda }$ but not in $H_{j,\rho }$.

\end{document}